\documentclass{article}
\usepackage{amssymb}
\usepackage{amsmath}
\usepackage{cite}
\usepackage{amsfonts}
\usepackage{amscd}
\usepackage[matrix]{xypic}
\newtheorem{theorem}{Theorem}

\newtheorem{defi}[theorem]{D\'efinition}
\newtheorem{theo}[theorem]{Th\'eor\`eme}
\newtheorem{coro}[theorem]{Corollaire}

\newtheorem{proposition}[theorem]{Proposition}

\newcommand{\g}{\frak{g}}

\newcommand{\z}{(\Bbb{Z}_2)^2}
\newcommand{\Z}{\Bbb{Z}}
\newcommand{\T}{\Bbb{T}}

\begin{document}

\title{\bf L'espace Riemannien et pseudo-Riemannien non sym\'etrique $SO(2m)/Sp(m)$.}

\author{Elisabeth REMM \thanks{%
corresponding author: e-mail: E.Remm@uha.fr} - Michel
GOZE \thanks{M.Goze@uha.fr.}\\
\\
{\small Universit\'{e} de Haute Alsace, F.S.T.}\\
{\small 4, rue des Fr\`{e}res Lumi\`{e}re - 68093 MULHOUSE - France}}
\date{}
\maketitle


\bigskip

\noindent {\bf Mots cl\'es.} Homogeneous manifolds, Riemannian structures, non symmetric spaces.

\begin{abstract}
In this work, we are interested in a non symmetric homogeneous space, namely $SO(2m)/Sp(m)$. 
We show that this space admits a structure of $\z$-symmetric space. We describe all the non degenerated
metrics and classify the Riemannian and Lorentzian ones.

\medskip

Dans ce travail, nous nous int\'eressons \`a un espace homog\`ene non sym\'etrique, \`a savoir $SO(2m)/Sp(m)$.  
Nous montrons 
que cet espace admet une structure d'espace $\z$-sym\'etrique. 
Nous d\'ecrivons toutes les m\'etriques non d\'eg\'en\'er\'ees et classons les m\'etriques riemanniennes et lorentziennes.

\end{abstract}

\section{Rappel. Espaces $\Gamma $-sym\'etriques}

Soit $\Gamma $ un groupe ab\'elien fini. Un espace $\Gamma $-sym\'etrique est un espace homog\`ene $M=G/H$
o\`u $G$ est un groupe de Lie connexe et $H$ un  sous-groupe ferm\'e de $G$ tel qu'il existe un homomorphisme
injectif de groupes
$$\rho: \gamma \in\Gamma \longrightarrow \rho (\gamma )\in  Aut(G),$$ 
tel que le groupe $H$ v\'erifie  $(G^\Gamma)_1 \subset H \subset G^\Gamma $ o\`u
$G^\Gamma =\left\{  g\in G / \forall \gamma \in \Gamma , \rho (\gamma )(g)=g \right\}$
et $(G^\Gamma )_1$ sa composante connexe passant par l'\'el\'ement neutre $1$ de $G$. On a en particulier:
$$
\left\{
\begin{array}{l}
\rho (\gamma _1)\circ \rho (\gamma _2)=\rho (\gamma _1 \cdot \gamma _2), \\
\\
 \rho (e_{\Gamma})=Id  \ \ {\mbox{\rm o\`u}} \ e_{\Gamma}  \ {\mbox{\rm est l'\'el\'ement neutre de }} \ \Gamma , \\
\\
\forall \gamma \in \Gamma, \rho (\gamma) (g)=g \   \Longleftrightarrow g\in H \ {\mbox{\rm d\`es que $H$ est connexe}} .
\end{array}
\right.
$$
Si $M=G/H$ est un espace $\Gamma $-sym\'etrique, alors l'alg\`ebre  de Lie $\frak{g}$ de $G$ est $\Gamma $-gradu\'ee:
$$\frak{g}=\oplus_{\gamma  \in \Gamma }  \frak{g}_\gamma $$
avec 
$$\left\{ 
\begin{array}{l}
[ \frak{g}_{\gamma_1}, \frak{g}_{\gamma_2}] \subset  \frak{g}_{\gamma_1 \cdot \gamma_2}  \\
  \frak{g}_{e_{\Gamma }}=\frak{h}
\end{array}
\right.$$
o\`u $\frak{h}$ est l'alg\`ebre de Lie de $H$. On dit que la paire $(\frak{g},\frak{h})$ est l'espace 
$\Gamma $-sym\'etrique local de $G/H$. Si $G$ est simplement connexe alors toute alg\`ebre de Lie $\Gamma $-gradu\'ee
$\frak{g}=\oplus_{\gamma  \in \Gamma }  \frak{g}_\gamma $ d\'efinit une paire $(\frak{g},\frak{g}_e)$ qui est l'espace 
$\Gamma $-sym\'etrique local d'un espace $\Gamma $-sym\'etrique $G/H$.

\medskip

Soit $M=G/H$ un espace $\Gamma $-sym\'etrique. En tout point $x$ de $M$ on peut d\'efinir un sous-groupe
$\Gamma _x$ de $\mathcal{D}iff(M)$ isomorphe \`a $\Gamma $ tel que $x$ soit le seul point fixe commun \`a tous les \'el\'ements 
$s_{\gamma ,x}$ de $\Gamma_x$. 

\smallskip

Les premiers exemples d'espaces $\Gamma$-sym\'etriques lorsque $\Gamma$ n'est pas cyclique correspondent au cas o\`u 
$\Gamma = \z$. Rappelons que tout espace $\Bbb{Z}_2$-sym\'etrique n'est rien d'autre qu'un espace sym\'etrique.

\medskip

\noindent {\bf Exemple. La sph\`ere $S^3$}.

Nous savons que la sph\`ere $S^3$ (ou plus g\'en\'eralement $S^n$) est munie d'une structure d'espace sym\'etrique. 
Il suffit de  consid\'erer un diff\'eomorphisme de $S^n$ sur l'espace homog\`ene $SO(n+1)/SO(n)$. 
L'alg\`ebre de Lie $so(n+1)$ admet une $\Z_2$-graduation
$$so(n+1)=so(n)\oplus \frak{m}$$
o\`u $so(n)$ est l'ensemble des points fixes de l'automorphisme involutif
$$
\tau _a: so(n+1)  \rightarrow   so(n+1) $$
donn\'e par
$$\tau _a(A)= SAS^{-1} $$
avec $S= \left(
\begin{array}{ll}
-1 & 0_n \\
0_n & I_n
\end{array} \right)$ et $I_n$ est la matrice identit\'e d'ordre $n$.

Nous pouvons munir la sph\`ere $S^3$ d'une structure d'espace $\z$- sym\'etrique en consid\'erant 
un diff\'eomorphisme de $S^3$ sur l'espace homog\`ene $SO(4)/Sp(2)$. Consid\'erons sur l'alg\`ebre de Lie $so(4)$
les automorphismes $\tau_a,\tau_b,\tau_c=\tau_a \circ \tau_b$ donn\'es par 
$$\left\{ 
\begin{array}{llll}
\tau_a:& M \in SO(4) & \mapsto J_a^{-1}MJ_a \\
\tau_b:& M \in SO(4) & \mapsto J_b^{-1}MJ_b \\
\end{array}
\right. $$
o\`u 
$$J_a=\left( 
\begin{array}{llll}
0&-1&0&0 \\
1&0&0&0 \\
0&0&0&1 \\
0&0&-1&0
\end{array}
\right) \ \mbox{\rm et} \  J_b=\left( 
\begin{array}{llll}
0&0&0&1 \\
0&0&-1&0 \\
0&1&0&0 \\
-1&0&0&0
\end{array}
\right) .$$ La famille $\left\{ Id,\tau_a,\tau_b,\tau_c \right\}$ d\'etermine un sous-groupe de $Aut(so(4))$ 
isomorphe \`a $\z$. Si 
$$\frak{g}_a=\left\{ M \in so(4) \, / \,  \tau_a(M)=M, \tau_b(M)=-M \right\},$$
$$\frak{g}_b=\left\{ M \in so(4) \, / \,  \tau_a(M)=-M, \tau_b(M)=M \right\}$$ et 
$$\frak{g}_c=\left\{ M \in so(4) \, / \,  \tau_a(M)=-M, \tau_b(M)=-M \right\},$$ on a la d\'ecomposition 
$\z$-sym\'etrique de $so(4)$:
$$so(4)=sp(2) \oplus \frak{g}_a  \oplus \frak{g}_b \oplus \frak{g}_c.$$
En effet l'ensemble de points fixes pour $\tau_a,\tau_b$ et $\tau_c$ est la sous-alg\`ebre de $so(4)$ dont
les \'el\'ements sont les matrices 
$$\left(
\begin{array}{rrrr}
0 & -a_2 & -a_3 & -a_4 \\
a_2 &  0& -a_4 &  a_3 \\
a_3 &  a_4 &  0 & -a_2 \\
a_4 & -a_3 & a_2 &  0 \\
\end{array}
\right).$$
Cette sous-alg\`ebre est isomorphe \`a $sp(2)$. Ainsi $(so(4),sp(2))$ est un espace local 
$\z$-sym\'etrique d\'eterminant une structure d'espace $\z$-sym\'etrique sur la shp\`ere $$S^3=SO(4)/Sp(2).$$

La structure d'espace sym\'etrique sur $S^3$ est li\'ee \`a l'existence en tout point $x \in S^3$
d'une sym\'etrie $s_x \in \mathcal{D}iff (S^3)$ v\'erifiant $s_x^2=Id$ et $x$ est le seul point 
tel que $s_x(x')=x'$ (seul point fixe).

La structure d'espace $\z$-sym\'etrique de $S^3$ donne l'existence d'un sous-groupe $\Gamma_x$ de $ \mathcal{D}iff(S^3) $
de "sym\'etries" de $S^3$:
$$\Gamma_x=\left\{ id, s_{a,x},s_{b,x} ,s_{c,x}  \right\}$$
o\`u $$\left\{ 
\begin{array}{l}
s_{x,a}^2=s_{x,b}^2=s_{x,c}^2 =Id \\
s_{x,a} \circ s_{x,b}=s_{x,b} \circ s_{x,a} =s_{x,c} \\
s_{x,b} \circ s_{x,c}=s_{x,c} \circ s_{x,b} =s_{x,a} \\
s_{x,a} \circ s_{x,c}=s_{x,c} \circ s_{x,a} =s_{x,b} 
\end{array}
\right.$$
et $x$ est le seul point fixe commun \`a toutes les sym\'etries:
$$\left( s_{x,a}(x')=s_{x,b}(x')=s_{x,c}(x')=x' \right) \Leftrightarrow x=x'.$$
Notons que chacune des sym\'etries peut avoir un ensemble de points fixes non r\'eduit \`a $\left\{ x \right\}$ et donc
chacune des sym\'etries ne d\'eterminent pas une structure d'espace sym\'etrique sur $S^3$. 
D\'eterminons ces sym\'etries.

Pour tout $\gamma \in \Gamma=\z=\left\{ e,a,b,c \right\}$, consid\'erons l'automorphisme de $SO(4)$ donn\'e par 
$$\left\{
\begin{array}{l}
\rho_a(A)=J_a^{-1}AJ_a \\ 
\rho_b(A)=J_b^{-1}AJ_b \\ 
\rho_c(A)=\rho_a \circ \rho_b(A)  \\ 
\rho_e(A)=A 
\end{array}
 \right. . $$
Si $x=[ A ]$ d\'esigne la classe dans l'espace homog\`ene $SO(4)/Sp(2)$ de la marice $A \in SO(4)$, alors 
$s_{a,x}[A]=[J_a^{-1}AJ_a], \ s_{b,x}[A]=[J_b^{-1}AJ_b ], \ s_{c,x}[A]=[J_c^{-1}AJ_c],$ o\`u $J_c=J_aJ_b$.
Ainsi chacune des sym\'etries a un grand cercle comme vari\'et\'e de points invariants.

\section{Espaces riemanniens $\Gamma$-sym\'etriques}

Soit $(M=G/H, \Gamma)$ un espace homog\`ene $\Gamma$-sym\'etrique.

\begin{defi}
Une m\'etrique riemannienne $g$ sur $M$ est dite adapt\'ee \`a la structure $\Gamma$-sym\'etrique si 
chacune des sym\'etries $s_{\gamma,x}$ est une isom\'etrie.
\end{defi}  

Si $\bigtriangledown_g $ est la connexion de Levi-Civita de $\frak{g}$, cette connexion ne co\"{\i}ncide pas en 
g\'en\'eral avec la connexion canonique $\bigtriangledown $ de l'espace homog\`ene ($\Gamma$-sym\'etrique).
Ces deux connexions co\"{\i}ncident si et seulement si $g$ est naturellement r\'eductive.

Par exemple dans le cas de la sph\`ere $S^3$ consid\'er\'ee comme espace $\z$-sym\'etrique, 
les m\'etriques adapt\'ees \`a cette structure sont les m\'etriques sur $SO(4)/Sp(2)$ invariantes par $SO(4)$
chacune \'etant
d\'efinie par une forme bilin\'eaire sym\'etrique $B$ sur $so(4)$ qui est $ad(sp(2))$-invariante. Si 
$so(4)=sp(2) \oplus \frak{g}_a \oplus \frak{g}_b \oplus  \frak{g}_c$ est 
la d\'ecomposition $\z$-gradu\'ee correspondante, le fait de dire que sur $S^3$ les sym\'etries 
$s_{\gamma,x}$ sont des isom\'etries est \'equivalent \`a dire que les espaces 
$\frak{g}_e,\frak{g}_a,\frak{g}_b, \frak{g}_c$ sont deux \`a deux orthogonaux pour $B$. D\'ecrivons en d\'etail cette graduation:
$$sp(2)=\left\{ \left(
\begin{array}{llll}
0 & -a_2 & -a_3 & -a_4 \\
a_2 &  0 & -a_4 &  a_3 \\
a_3 &  a_4 &  0 & -a_2 \\
a_4 & -a_3 & a_2 &  0 \\
\end{array}
\right)
\right\}, \ \frak{g}_a=\left\{ \left(
\begin{array}{llll}
0 & 0 & 0 & x \\
0&  0 & -x & 0 \\
0 &  x &  0 & 0 \\
-x & 0 & 0 &  0 \\
\end{array}
\right)
\right\},$$
$$\frak{g}_b=\left\{ \left(
\begin{array}{llll}
0 & y & 0 & 0 \\
-y &  0 & 0 & 0 \\
0 &  0 &  0 & -y \\
0 & 0 & y &  0 \\
\end{array}
\right)
\right\} \ \mbox{\rm et} \ \frak{g}_c=\left\{ \left(
\begin{array}{llll}
0 & 0 & z & 0 \\
0&  0 & 0 & z \\
-z &  0 &  0 & 0 \\
0 & -z & 0 &  0 \\
\end{array}
\right)
\right\}.$$ 
Si $\left\{ A_1 ,A_2, A_3, X,Y,Z \right\}$ est une base adapt\'ee \`a cette graduation et si 
$\left\{ \alpha _1 ,\alpha _2, \alpha _3, \omega_1, \omega_2, \omega_3 \right\}$ en est la base duale alors
$$B\mid_{ \frak{g}_a \oplus \frak{g}_b \oplus \frak{g}_c  } = 
\lambda _1^2\omega _1^2+\lambda _2^2\omega _2^2+\lambda _3^2\omega _3^2.$$
La m\'etrique correspondante sera naturellement r\'eductive si et seulement si $\lambda_1= \lambda_2 = \lambda_3$
et dans ce cas-l\`a elle correspond \`a la restriction de la forme de Killing Cartan.

\section{M\'etriques adapt\'ees \`a la structure $\z$-sym\'etrique de $SO(2m)/Sp(m)$}

\subsection{La graduation $\z$-sym\'etrique}

Consid\'erons les matrices 
$$S_m= \left(
\begin{array}{ll}
0& I_m  \\
-I_n & 0
\end{array} \right),  \ X_a=\left(
\begin{array}{ll}
-1& 0  \\
 0 & 1
\end{array} \right),X_b=\left(
\begin{array}{ll}
 0&1  \\
 1 & 0
\end{array} \right),X_c=\left(
\begin{array}{ll}
0& -1  \\
 1 & 0
\end{array} \right).$$
Soit $M \in so(2m).$ Les applications 
$$\begin{array}{l}
\tau_a(M)=J_a^{-1}MJ_a   \\
\tau_b(M)=J_b^{-1}MJ_b  \\
\tau_a(M)=J_c^{-1}MJ_c  
\end{array}$$
o\`u $J_a=S_m \otimes X_a, \ J_b=S_m \otimes X_b, \ J_c=S_m \otimes X_c$ sont des automorphismes involutifs de $so(2m)$
qui commutent deux \`a deux. Ainsi $\left\{ Id,\tau _a,\tau _b,\tau _c \right\}$ est un sous groupe de 
$Aut (so(2m))$ isomorphe \`a $\z$. Il d\'efinit donc une $\z$-graduation 
$$so(2m)=   \frak{g}_e \oplus  \frak{g}_a \oplus \frak{g}_b \oplus \frak{g}_c$$
o\`u  
$$\begin{array}{l}
\frak{g}_e= \left\{ M \in so(2m)\, / \, \tau _a(M)=\tau _b(M)=\tau _c (M)=M \right\} \\ 
\frak{g}_a= \left\{ M \in so(2m)\, / \, \tau _a(M) =\tau _c (M)=-M, \tau _b(M)=M \right\} \\ 
\frak{g}_b= \left\{ M \in so(2m)\, / \, \tau _b(M)=\tau _c(M)=-M, \tau _a (M)=M \right\} \\ 
\frak{g}_c= \left\{ M \in so(2m)\, / \, \tau _a(M)=\tau _b(M)=-M, \tau _c (M)=M \right\}  
\end{array}$$
Ainsi
$$
\frak{g}_e= \left\{ \left(
\begin{array}{rr|rr}
A_1 & B_1 & A_2 & B_2 \\
-B_1& A_1 & B_2 & -A_2 \\
\hline 
-^t \! A_2 & -^t \! B_2 & A_1 & B_1 \\
-^t \! B_2 & ^t \! A_2 & -B_1 & A_1 
\end{array}
\right)
\ \mbox{\rm avec} \ 
\begin{array}{cc}
^t \! A_1 =  -A_1, & ^t \! B_1  =  B_1 \\
^t \! A_2  =  A_2, & ^t \! B_2  =  B_2 \\
\end{array}
\right\}$$
$$
\frak{g}_a= \left\{ \left(
\begin{array}{rr|rr}
X_1 & Y_1 & Z_1 & T_1 \\
Y_1& -X_1 & -T_1 & Z_1 \\
\hline 
-^t \! Z_1 & ^t \! T_1 & -X_1 & -Y_1 \\
-^t \! T_1 & -^t \! Z_1 & -Y_1 & X_1 
\end{array}
\right)
\ \mbox{\rm avec} \ 
\begin{array}{cc}
^t \! X_1 =  -X_1, & ^t \! Y_1  =  -Y_1 \\
^t \! Z_1  =  -Z_1, & ^t \! T_1  =  T_1 \\
\end{array}
\right\} \label{g_a}$$
$$
\frak{g}_b= \left\{ \left(
\begin{array}{rr|rr}
X_2 & Y_2 & Z_2 & T_2 \\
-Y_2& X_2 & T_2 & Z_2 \\
\hline 
-^t \! Z_2 & -^t \! T_2 & -X_2 & -Y_2 \\
-^t \! T_2 & -^t \! Z_2 & Y_2 & -X_2 
\end{array}
\right)
\ \mbox{\rm avec} \ 
\begin{array}{cc}
^t \! X_2 =  -X_2, & ^t \! Y_2  =  Y_2 \\
^t \! Z_2  =  -Z_2, & ^t \! T_2  =  -T_2 \\
\end{array}
\right\}$$
$$
\frak{g}_c= \left\{ \left(
\begin{array}{rr|rr}
X_3 & Y_3 & Z_3 & T_3 \\
Y_3& -X_3 & -T_3 & Z_3 \\
\hline 
-^t \! Z_3 & ^t \! T_3 & X_3 & Y_3 \\
-^t \! T_3 & -^t \! Z_3 & Y_3 & -X_3 
\end{array}
\right)
\ \mbox{\rm avec} \ 
\begin{array}{cc}
^t \! X_3=  -X_3, & ^t \! Y_3  =  -Y_3 \\
^t \! Z_3  =  Z_3, & ^t \! T_3  =  -T_3 \\
\end{array}
\right\}$$
Notons que $dim \frak{g}_e=m(2m+1), \, dim \frak{g}_a=dim \frak{g}_b=dim \frak{g}_c=m(2m-1)$.

\begin{proposition}
Dans cette graduation $\frak{g}_e$ est isomorphe \`a $sp(m)$ et toute $\z$-graduation de $so(2m)$ telle que 
$\frak{g}_e$ soit isomorphe \`a $sp(m)$ est \'equivalente \`a la graduation ci-dessus.
\end{proposition}

En effet $\frak{g}_e$ est simple de rang $m$ et de dimension $m(2m+1)$. 
La deuxi\`eme partie r\'esulte de la classification donn\'ee dans \cite{B.G} et \cite{B.G.R}. 

\begin{coro}
Il n'existe, \`a \'equivalence pr\`es, qu'une seule structure d'espace homog\`ene $\z$-sym\'etrique sur 
l'espace homog\`ene compact $SO(2m)/Sp(m)$.
\end{coro}

Cette structure est associ\'ee \`a l'existence en tout point $x$ de $SO(2m)/Sp(m)$ 
d'un sous-groupe de $\mathcal{D}iff(M)$ isomorphe \`a $\z$. Notons $\Gamma_x$ ce sous-groupe. Il est enti\`erement
d\'efini d\`es que l'on connait $\Gamma_{\bar{1} }$ o\`u $\bar{1}$ est la classe dans $SO(2m)/Sp(m)$ de l'\'el\'ement
neutre $1$ de $SO(2m)$. Notons 
$$\Gamma_{\bar{1} }=\left\{ s_{e,\bar{1}}, s_{a,\bar{1}}, s_{b,\bar{1}}, s_{c,\bar{1}} \right\},$$ 
les sym\'etries $s_{\gamma ,\bar{1}}(x)=\pi(\rho_\gamma(A))$ o\`u $\pi:SO(2m) \rightarrow SO(2m)/Sp(m)$
est la submersion canonique, $x=\pi(A)$ et $\rho_\gamma$ est un automorphisme de $SO(2m)$ dont l'application tangente en $1$ 
coincide avec $\tau_\gamma$. Ainsi
$$\left\{
\begin{array}{l}
\rho _a(A)=J_a^{-1}AJ_a \\  
\rho _b(A)=J_b^{-1}AJ_b \\  
\rho _c(A)=J_c^{-1}AJ_c \\  
\end{array}
\right. .$$
Si $B \in \pi(A)$ alors il existe $P \in Sp(m)$ tel que $B=AP$. 
On a $J_a^{-1}BJ_a = J_a^{-1}AJ_aJ_a^{-1}PJ_a = J_a^{-1}AJ_a$ car $P$ est invariante pour tous les automorphismes 
$\rho_a,\rho _b,\rho _c.$

\subsection{Structure m\'etrique $\z$-sym\'etrique }

\noindent Une m\'etrique non d\'eg\'en\'er\'ee $g$ invariante par $SO(2m)$ sur $SO(2m)/Sp(m)$ est  adapt\'ee \`a la $\z$-structure
si les sym\'etries $s_{x,\gamma}$ sont des isom\'etries c'est-\`a-dire si les automorphismes $\rho _\gamma$
induisent des isom\'etries lin\'eaires.

 \noindent Ceci implique que $g$ soit d\'efinie par une forme bilin\'eaire sym\'etrique non d\'eg\'en\'er\'ee
$B$ sur $\frak{g}_a \oplus \frak{g}_b \oplus \frak{g}_c $ telle que les espaces 
$\frak{g}_a,\, \frak{g}_b,\, \frak{g}_c$ soient deux \`a deux orthogonaux. D\'eterminons toutes les formes 
bilin\'eaires $B$ v\'erifiant les hypoth\`eses ci-dessus. Une telle forme s'\'ecrit donc
$$B=B _a +B _b+B _c$$
o\`u $B _a$(resp. $B _b$, resp. $B _c$ ) est une forme bilin\'eaire sym\'etrique non d\'eg\'en\'er\'ee invariante par 
$\frak{g}_e$ dont le noyau contient $\frak{g}_b \oplus \frak{g}_c$ 
 (resp. $\frak{g}_a \oplus \frak{g}_c$, resp. $\frak{g}_a \oplus \frak{g}_b$).  
   
\subsection{Exemples}

1) Dans le cas de la sph\`ere $SO(4)/Sp(2)$ la m\'etrique adapt\'ee \`a la structure $\z$-sym\'etrique est d\'efinie
par la forme bilin\'eaire $B$ sur $\frak{g}_a \oplus \frak{g}_b \oplus \frak{g}_c$ qui est 
$ad(sp(2))$-invariante. Nous avons vu qu'une telle forme s'\'ecrivait
$$ B=\lambda _1^2\omega _1^2+\lambda _2^2\omega _2^2+\lambda _3^2\omega _3^2.$$
Elle est d\'efinie positive si et seulement si les coefficients $\lambda_i$ sont positifs ou nuls.
  
\medskip
 
2) Consid\'erons l'espace $\z$-sym\'etrique compact $SO(8)/Sp(4)$. Afin de fixer les notations \'ecrivons la
$\z$-graduation de $so(8)$ ainsi:
$$
\frak{g}_a= \left\{ \left(
\begin{array}{rr|rr}
X_1 & Y_1 & Z_1 & T_1 \\
Y_1& -X_1 & -T_1 & Z_1 \\
\hline 
-^t \! Z_1 & ^t \! T_1 & -X_1 & -Y_1 \\
-^t \! T_1 & -^t \! Z_1 & -Y_1 & X_1 
\end{array}
\right)
\ \mbox{\rm avec} \ 
\begin{array}{cc}
^t \! X_1 =  -X_1, & ^t \! Y_1  =  -Y_1 \\
^t \! Z_1  =  -Z_1, & ^t \! T_1  =  T_1 \\
\end{array}
\right\}$$
$$
\frak{g}_b= \left\{ \left(
\begin{array}{rr|rr}
X_2 & Y_2 & Z_2 & T_2 \\
-Y_2& X_2 & T_2 & Z_2 \\
\hline 
-^t \! Z_2 & -^t \! T_2 & -X_2 & -Y_2 \\
-^t \! T_2 & ^t \! Z_2 & Y_2 & -X_2 
\end{array}
\right)
\ \mbox{\rm avec} \ 
\begin{array}{cc}
^t \! X_2 =  -X_2, & ^t \! Y_2  =  Y_2 \\
^t \! Z_2  =  -Z_2, & ^t \! T_2  =  -T_2 \\
\end{array}
\right\}$$
$$
\frak{g}_c= \left\{ \left(
\begin{array}{rr|rr}
X_3 & Y_3 & Z_3 & T_3 \\
Y_3& -X_3 & -T_3 & Z_3 \\
\hline 
-^t \! Z_3 & ^t \! T_3 & X_3 & Y_3 \\
-^t \! T_3 & -^t \! Z_3 & Y_3 & -X_3 
\end{array}
\right)
\ \mbox{\rm avec} \ 
\begin{array}{cc}
^t \! X_3=  -X_3, & ^t \! Y_3  =  -Y_3 \\
^t \! Z_3  =  Z_3, & ^t \! T_3  =  -T_3 \\
\end{array}
\right\}$$
et pour la matrice $X_i$ (resp. $Y_i,Z_i,T_i$) on notera 
$X_i=\left(
\begin{array}{cc}
0&x_i \\
-x_i & 0
\end{array}
\right)$ si elle est antisym\'etrique 
ou $X_i=\left(
\begin{array}{cc}
x_i^1 &x_i^2 \\
x_i^2 & x_i^3
\end{array}
\right)$ si elle est sym\'etrique, c'est \`a dire $X_i=\sum_j x_i^jX_i^j .$

Enfin on notera par les lettres $\alpha_i,\beta_i, \gamma_i, \delta_i$ les formes lin\'eaires duales des 
vecteurs d\'efinis respectivement par les matrices $X_i,Y_i,Z_i,T_i$. Ainsi si $X_i$ est antisym\'etrique, 
la forme duale correspondante sera not\'ee $\alpha_i$, et si $X_i$ est sym\'etrique, 
les formes duales $ \alpha_i^1,\alpha_i^2,\alpha_i^3$ correspondent aux vecteurs 
$\left(
\begin{array}{cc}
1&0 \\
0&0
\end{array}
 \right),\left(
\begin{array}{cc}
0&1 \\
1&0
\end{array}
 \right),\left(
\begin{array}{cc}
0&0 \\
0&1
\end{array}
 \right) .$ 
Ceci \'etant la forme 
$B$ s'\'ecrit $B_a+B_b+B_c$ o\`u la forme 
$B_{\gamma}$ a pour noyau
$\frak{g}_{\gamma_1} \oplus \frak{g}_{\gamma_2}$ avec $\frak{g}_{\gamma} \neq \frak{g}_{\gamma_1}$ et 
$\frak{g}_{\gamma} \neq \frak{g}_{\gamma_2}$. D\'eterminons $B_a$. Comme elle est invariante par $ad(sp(2))$
on obtient:
$$\left\{
\begin{array}{l}
\medskip
B_a(X_1,Y_1)=B_a(X_1,Z_1)=B_a(X_1,T_1^i)=0 \\
\medskip
B_a(Y_1,Z_1)=B_a(Y_1,T_1^i)=0 \\
\medskip
B_a(Z_1,T_1^i)=B_a(T_1^i,T_1^2)=0 \mbox{\rm pour } \, i=1,3 \\
\medskip
B_a(X_1,X_1)=B_a(Y_1,Y_1)=B_a(Z_1,Z_1)=B_a(T_1^2,T_1^2)=0 \\
\medskip
B_a(T_1^1,T_1^1)=B_a(T_1^3,T_1^3) \\
\medskip
B_a(X_1,X_1)=2B_a(T_1^1,T_1^1)-2B_a(T_1^1,T_1^3)
\end{array}
\right.$$
Ainsi la forme quadratique associ\'ee s'\'ecrit
$$q_{\g_a}=\lambda _1(\alpha _1^2+\beta _1^2+\gamma _1^2+(\delta _1^2)^2)+
\lambda _2((\delta _1^1)^2)+(\delta _1^3)^2)+(\lambda_2-\frac{ \lambda_1}{2} )((\delta _1^1)(\delta _1^3))$$
soit
$$q_{\g_a}=\lambda _1(\alpha _1^2+\beta _1^2+\gamma _1^2+(\delta _1^2)^2)+
(\frac{3\lambda _2}{4}-\frac{\lambda _1}{8})(\delta _1^1+\delta _1^3)^2+
(\frac{\lambda _2}{4}+\frac{\lambda _1}{8})(\delta _1^1-\delta _1^3)^2.$$
De m\^eme nous aurons
$$q_{\g_b}=\lambda _3(\alpha _2^2+(\beta _2^2)^2+\gamma _2^2+\delta _2^2)+
(\frac{3\lambda _4}{4}-\frac{\lambda _3}{8})(\beta _2^1+\beta_2^3)^2+
(\frac{\lambda _4}{4}+\frac{\lambda _3}{8})(\beta_2^1-\beta_2^3)^2$$
et
$$q_{\g_c}=\lambda _5(\alpha _3^2+\beta _3^2+(\gamma _3^2)^2+\delta _3^2)+
(\frac{3\lambda _6}{4}-\frac{\lambda _5}{8})(\gamma _3^1+\gamma_3^3)^2+
(\frac{\lambda _6}{4}+\frac{\lambda _5}{8})(\gamma_3^1-\gamma_3^3)^2.$$

\noindent {\bf Remarques.}
1. La forme $B$ d\'efinit une m\'etrique riemannienne si et seulement si 
$$\lambda_{2p}> \frac{\lambda _{2p-1}}{6}>0$$
pour $p=1,2,3$.
Si cette  contrainte est relach\'ee, la forme $B$, suppos\'ee non d\'eg\'en\'er\'ee, peut d\'efinir une m\'etrique pseudo riemannienne
sur l'espace $\z$-sym\'etrique. Nous verrons cela dans le dernier paragraphe.

2. Consid\'erons le sous-espace $\frak{g}_e \oplus \frak{g}_a$. Comme $[\frak{g}_a,\frak{g}_a] \subset \frak{g}_e$
c'est un sous-alg\`ebre de $so(8)$ (ou plus g\'en\'eralement de $\frak{g}$) admettant une stucture sym\'etrique.
La forme $B_a$ induit donc une structure riemannienne ou pseudo-riemannienne sur l'espace sym\'etrique associ\'e
\`a l'espace sym\'etrique local 
$(\frak{g}_e , \frak{g}_a)$.  Dans l'exemple pr\'ec\'edent 
$\frak{g}_e \oplus \frak{g}_a$ est la sous-alg\`ebre de $so(8)$ donn\'ee par les matrices: 
$$
\left(
\begin{array}{rr|rr}
X_1 & X_3 & X_4 & X_5 \\
-^t \! X_3& X_2 & X_6 & -^t \! X_4 \\
\hline 
-^t \! X_4 & -^t \! X_6 & X_2 & ^t \! X_3 \\
-^t \! X_5 &  X_4 & -^t \! X_3 & X_2 
\end{array}
\right)
\ \mbox{\rm avec} \ 
\begin{array}{cc}
^t \! X_1 =  -X_1, & ^t \! X_2  =  -X_2 \\
^t \! X_5  =  X_5, & ^t \! X_6  =  X_6 \\
\end{array}
$$
Dans \cite{B}, on d\'etermine les espaces r\'eels en \'etudiant ces structures sym\'etriques 
$\frak{g}_e \oplus \frak{g}_a$ donn\'ees par deux automorphismes commutant de $\frak{g}$. En effet
si $\frak{g}$ est simple r\'eelle et si $\sigma $ est un automorphisme involutif de $ \frak{g} $, 
il existe une sous-alg\`ebre compacte maximale $\frak{g}_1$  de $\frak{g}$ qui est invariante 
par $\sigma $ et l'\'etude des espaces locaux sym\'etriques $(\frak{g} ,\frak{g}_e)$
se ram\`ene \`a l'\'etude des espaces locaux sym\'etriques $(\frak{g}_1 ,\frak{g}_{11})$ o\`u $\frak{g}_1$ est 
compacte. Dans ce cas  $ \frak{g} $ est d\'efinie \`a partir de $ \frak{g}_1$ par un automorphisme involutif $\tau$
commutant avec l'automorphisme $\sigma$.
Ici notre approche est en partie similaire mais le but est de regarder la structure des espaces non sym\'etrique 
associ\'es aux paires $( \frak{g},\frak{g}_e ).$ 

 \smallskip

Dans le cas particulier de l'espace $\z$-sym\'etrique compact $SO(8)/Sp(4)$ l'alg\`ebre de Lie 
$ \frak{g}_e \oplus  \frak{g}_a $ est isomorphe \`a $so(4) \oplus \mathbb{R}$ o\`u $\mathbb{R}$ d\'esigne 
l'alg\`ebre ab\'elienne de dimension $1$. Notons \'egalement que chacun des espaces sym\'etriques
$ \frak{g}_e \oplus  \frak{g}_a, \frak{g}_e \oplus  \frak{g}_b, \frak{g}_e \oplus  \frak{g}_c$ est isomorphe
\`a $so(4) \oplus \mathbb{R}$. Mais ceci n'est pas g\'en\'eral, les alg\`ebres sym\'etriques peuvent ne pas \^etre
isomorphes  ni m\^eme de m\^eme dimension. L'espace sym\'etrique compact connexe associ\'e est l'espace homog\`ene 
${Su(4)}/{Sp(2)} \times \T$ o\`u $\T$ est le tore \`a une dimension. C'est un espace riemannien sym\'etrique compact non irr\'eductible. La m\'etrique
$q_{\frak{g}_a}$ d\'efinie pr\'ec\'edemment correspond \`a une m\'etrique riemannienne ou pseudo riemannienne sur cet 
espace. La restriction au premier facteur correspond \`a la m\'etrique associ\'ee \`a la forme de Killing Cartan sur 
$su(4)$. Elle correspond \`a $\lambda_2=\frac{ \lambda_1}{2}$. 

\subsection{Cas g\'en\'eral: m\'etriques adapt\'ees sur $SO(2m)/Sp(m)$}

Notations. Nous avons \'ecrit une matrice g\'en\'erale de $\frak{g}_a$ sous la forme (\ref{g_a}).
Si on note $(X_1,Y_1,Z_1,T_1)$ un \'el\'ement de $\g_a$, on consid\`ere  la base de $\frak{g}_a$,
$\{X_{1,ij},Y_{1,ij},Z_{1,ij},T_{1,ij}\}$
correspondant aux matrices \'el\'ementaires . La base duale sera not\'ee 
$(\alpha _{a,ij},\beta_{a,ij}, \gamma_{a,ij}, \delta_{a,ij} )$. Rappelons que $X_1,Y_1,Z_1$ sont antisym\'etriques 
alors que $T_1$ est sym\'etrique. Les crochets correspondent aux repr\'esentations de $so(\frac{m}{2})$ sur 
lui-m\^eme ou de $so(\frac{m}{2})$ sur l'espace des matrices sym\'etriques. On aura donc 
$$q_{\frak{g}_a}=\lambda _1^a(\sum (\alpha _{a,ij}^2+\beta _{a,ij}^2+\gamma _{a,ij}^2)+\sum_{i \neq j}\delta _{a,ij}^2)+
\lambda _2^a(\delta _{a,ii}^2)+(\lambda_2^a-\frac{ \lambda_1^a}{2} )(\sum_{i<j}(\delta _{a,ii}\delta _{a,jj}).$$
Les formes $q_{\frak{g}_b}$ et $q_{\frak{g}_c}$ admettent une d\'ecomposition analogue, en tenant compte du fait que dans $\g_b$
ce sont les matrices $Y_2$ qui sont sym\'etriques et pour $\g_a$ les matrices $Z_1$ (\ref{g_a}). On note
$(\alpha _{b,ij},\beta_{b,ij}, \gamma_{b,ij}, \delta_{b,ij} )$ la base duale de $\{X_{2,ij},Y_{2,ij},Z_{2,ij},T_{2,ij}\}$
et par $(\alpha _{c,ij},\beta_{c,ij}, \gamma_{c,ij}, \delta_{c,ij} )$ la base duale de $\{X_{3,ij},Y_{3,ij},Z_{3,ij},
T_{3,ij}\}$.

\begin{proposition}
Toute m\'etrique non d\'eg\'en\'er\'e adapt\'ee \`a la structure $\z$-sym\'etrique de l'espace homog\`ene
$SO(2m)/Sp(m)$ est d\'efinie \`a partir de la forme bilin\'eaire $ad(\g_e)$-invariante sur $\g_a\oplus \g_b\oplus \g_c$
$B=q_{\frak{g}_a}+q_{\frak{g}_b}+q_{\frak{g}_b}$ avec
$$
\left\{
\begin{array}{l}
\medskip

q_{\frak{g}_a}=\lambda _1^a\left(\sum (\alpha _{a,ij}^2+\beta _{a,ij}^2+\gamma _{a,ij}^2\right)+\sum_{i \neq j}\delta _{a,ij}^2)+
\lambda _2^a(\delta _{a,ii}^2)+(\lambda_2^a-\frac{ \lambda_1^a}{2} )(\sum_{i<j}(\delta _{a,ii}\delta _{a,jj})\\

\medskip

q_{\frak{g}_b}=\lambda _1^b(\sum (\alpha _{b,ij}^2 +\gamma _{ij}^2)+\delta _{b,ij}^2+\sum_{i \neq j}\beta _{b,ij}^2)+
\lambda _2^b(\beta _{b,ii}^2)+(\lambda_2^b-\frac{ \lambda_1^b}{2} )(\sum_{i<j}(\beta _{b,ii}\beta _{b,jj})\\
\medskip
q_{\frak{g}_c}=\lambda _1^c(\sum (  \beta _{c,ij}^2+\gamma _{c,ij}^2)+\delta _{c,ij}^2+\sum_{i \neq j}\alpha _{c,ij}^2)+
\lambda _2^c(\alpha _{c,ii}^2)+(\lambda_2^c-\frac{ \lambda_1^c}{2} )(\sum_{i<j}(\alpha _{c,ii}\alpha _{c,jj})
\end{array}
\right.
$$
\end{proposition}

\section{M\'etriques  pseudo-riemanniennes $\z$-sym\'etriques sur $SO(2m)/Sp(m)$ } 

\subsection{Signature des formes $q_{\g{_\gamma} }$}
Soit $\gamma \in \{a,b,c\}$. Les valeurs propres de la forme $q_{\g_\gamma }$ sont
$$\mu_{1,\gamma }=\lambda_1^\gamma , \ \ \mu _{2,\gamma }=\lambda _2^\gamma /2+\lambda _{1}^\gamma /4, \ \ 
\mu_{3,\gamma }=\lambda _2^\gamma\frac{r+1}{2}-\lambda _{1}^\gamma \frac{r-1}{4}$$
o\`u $r$ est l'ordre commun des matrices sym\'etriques $X_4,Y_2,Z_1$. Ces valeurs propres sont respectivement de
multiplicit\'e $dim \g_{\gamma }-r,r-1,1$.
Le signe des valeurs propres $ \mu _{2,\gamma }$ et $\mu _{3,\gamma }$ est donc
$$
\left\{
\begin{array}{l}
 \mu _{2,\gamma } > 0 \Longleftrightarrow \lambda_2^\gamma >-\lambda_1^\gamma /2\\
\mu _{3,\gamma } > 0 \Longleftrightarrow  \lambda_2^\gamma > -\lambda_1^\gamma \frac{r-1}{2(r+1)}
\end{array}
\right.
$$
On en d\'eduit, si $s(q)$ d\'esigne la signature de la forme quadratique $q$ :
$$
\left\{
\begin{array}{lll}
\medskip
s(q_{\g_\gamma })& = (dim \g_\gamma ,0) & \Leftrightarrow ( \lambda_1^\gamma >0, 
\ \lambda_2^\gamma > \lambda_1^\gamma \frac{r-1}{2(r+1)})\\
\medskip
& = (dim \g_\gamma -1 ,1) & \Leftrightarrow ( \lambda_1^\gamma >0, \ -\lambda_1^\gamma /2< \lambda_2^\gamma 
<\lambda_1^\gamma \frac{r-1}{2(r+1)})\\
\medskip
& = (dim \g_\gamma -r ,r) & \Leftrightarrow ( \lambda_1^\gamma >0, \ \lambda_2^\gamma <-\lambda_1^\gamma /2)\\
\medskip
& = (r,dim \g_\gamma -r) & \Leftrightarrow ( \lambda_1^\gamma <0,  \ \lambda_2^\gamma >-\lambda_1^\gamma /2)\\
\medskip
& = (1,dim \g_\gamma -1) & \Leftrightarrow ( \lambda_1^\gamma <0,  \ \lambda_1^\gamma \frac{r-1}{2(r+1)}
< \lambda_2^\gamma <-\lambda_1^\gamma /2)\\
\medskip
& = (0, dim \g_\gamma ) & \Leftrightarrow ( \lambda_1^\gamma <0, \ \lambda_2^\gamma < \lambda_1^\gamma \frac{r-1}{2(r+1)}\\
\end{array}
\right.
$$ 
Notons que $\mu_{2,\gamma }= \mu_{3,\gamma }$ si et seulement si $\lambda_1^\gamma= 2\lambda_2^\gamma.$
\subsection{Classification des m\'etriques riemanniennes adapt\'ees sur $SO(2m)/Sp(m)$}
Comme $r=\frac{m^2+m-2}{m^2+m+2}$ on a le r\'esultat suivant
\begin{theo}
Toute m\'etrique riemannienne sur $SO(2m)/Sp(m)$ adapt\'ee \`a la structure $\z$-sym\'etrique
est d\'efinie
 \`a partir de la forme bilin\'eaire sur $\g_a\o^plus \g_b\oplus \g_c$ 
$$B= q_{\g_a }( \lambda_1^a,\lambda_2^a)+q_{\g_b }( \lambda_1^b,\lambda_2^b)+q_{\g_b }( \lambda_1^b,\lambda_2^b)$$
avec
$$\left\{
\begin{array}{l}
\lambda_1^\gamma >0\\
\lambda_2^\gamma > \lambda_1^\gamma \frac{m^2+m-2}{2(m^2+m+2})\\
\end{array}
\right.
$$
pour tout $\gamma \in \{a,b,c\}$.
\end{theo}
Pour une telle m\'etrique, la connexion de Levi-Civita ne co\"{\i}ncide pas en g\'en\'eral avec la connexion
canonique associ\'ee \`a la structure $\z$-sym\'etrique (\cite{B.G.R}). Ces deux connexions 
sont les m\^emes si et seulement si la m\'trique riemannienne est naturellement r\'eductive. Elle 
correspond donc \`a la restriction de la forme
de Killing (au signe pr\`es) de $SO(2m)$.  Cette m\'etrique correspond \`a la forme bilin\'eaire $B$ d\'efinie par les param\`etres
$$ \lambda_1^a= \lambda_1^b=\lambda_1^c= 2\lambda_2^a=2\lambda_2^b=2\lambda_2^c.$$

\subsection{Classification des m\'etriques lorentziennes adapt\'ees sur l'espace $SO(2m)/Sp(m)$}
Les m\'etriques lorentziennes adapt\'ees \`a la structure $\z$-sym\'etrique sont d\'efinies par les formes bilin\'eaires $B$, d\'efinies dans la section pr\'ec\'edente,
dont la signature est $(dim(\g_a)+dim(\g_b)+dim(\g_c)-1,1)$. On a donc

\begin{theo}
Toute m\'etrique lorentzienne sur $SO(2m)/Sp(m)$ adapt\'ee \`a la structure $\z$-sym\'etrique est d\'efinie par l'une des
 formes bilin\'eaires 
$$B= q_{\g_a }( \lambda_1^a,\lambda_2^a)+q_{\g_b }( \lambda_1^b,\lambda_2^b)+q_{\g_b }( \lambda_1^b,\lambda_2^b)$$
avec
$$
\left\{
\begin{array}{l}
\medskip
\forall \gamma \in \{a,b,c\}, \ \lambda_1^\gamma >0 \\
\medskip
\exists \gamma _0 \in \{a,b,c\} \ \mbox{\rm tel \ que} \ 
-\lambda_1^{\gamma_0} /2< \lambda_2^{\gamma_0} 
<\lambda_1^{\gamma_0} \frac{r-1}{2(r+1)}\\
\medskip
\forall \gamma  \neq \gamma _0, \ \ \lambda_2^\gamma > \lambda_1^\gamma \frac{r-1}{2(r+1)}. \ \\
\end{array}
\right.
$$
\end{theo}

\end{document}